\documentclass[12pt]{article}
\usepackage{amsmath}
\usepackage{amsfonts}
\usepackage{mathrsfs}
\usepackage{amssymb}
\usepackage[usenames]{color}
\usepackage{CJK}

\def\vspb{\vs{5pt}}

\def\qed{\hfill\mbox{$\Box$}}
\def\bs{\backslash}
\def\cl{\centerline}
\def\dis{\displaystyle}

\def\ni{\noindent}

\def\ol{\overline}
\def\ra{\hookrightarrow}
\def\sc{\scriptstyle}
\def\ssc{\scriptscriptstyle}

\def\vs{\vspace*}

\def\wt{\widetilde}

\numberwithin{equation}{section}
\newtheorem{theo}{Theorem}[section]
\newtheorem{defi}[theo]{Definition}

\newtheorem{lemm}[theo]{Lemma}
\newtheorem{ques}[theo]{Question}
\newtheorem{prop}[theo]{Proposition}
\newtheorem{rema}[theo]{Remark}

\newtheorem{clai}{Claim}
\newtheorem{conj}[theo]{Conjecture}

\def\a{\alpha}
\def\b{\beta}
\def\d{\delta}
\def\D{\Delta}

\def\g{\gamma}
\def\G{\Gamma}

\def\L{\Lambda}
\def\o{\omega}
\def\p{\partial}

\def\T{\Theta}

\def\Z{\mathbb{Z}}
\def\C{\mathbb{C}}
\def\Q{\mathbb{Q}}

\def\mub{{\mu^+}}
\def\WW{{\cal W}}
\def\BB{{\cal B}(q)}
\def\LL{{\cal L}(q,\mu_0)}

\def\PP{{\cal P}(q)}
\def\BBB{{\cal B}}
\def\LLL{{\cal L}}
\def\PPP{{\cal P}}
\def\sZ{{\ssc\,}\mathbb{Z}}
\def\Vir{{\rm{Vir}}}
\def\sp{{\rm{span}}}
\def\dim{{\rm{dim}}}
\def\max{{\rm{max}}}
\def\deg{{\rm{deg}}}
\def\LHS{{\theta_{\rm{lhs}}}}
\def\RHS{{\theta_{\rm{rhs}}}}
\def\Eq{{\rm{E\ssc}}}
\def\Eqa{{\rm{Ea\ssc}}}
\def\Eqb{{\rm{Eb\ssc}}}
\def\Eqc{{\rm{Ec\ssc}}}

\begin{document}
\begin{CJK*}{GBK}{song}
\cl{\large\bf {Quasifinite
representations 
of a class of Block type Lie algebras $\BB$}
\footnote{Supported by NSF grant 10825101 of China\\\indent\ \
Corresponding author: C.~Xia (chgxia@mail.ustc.edu.cn)}} \vskip5pt

\cl{{Yucai Su$^{\,\dag,\,\ddag}$, \ \ Chunguang Xia$^{\,\ddag}$, \ \ Ying Xu$^{\,\ddag}$}}

\cl{\small $^{\dag\,}$Department of Mathematics, Tongji
University, Shanghai 200092, \vs{2pt}China} \cl{\small
$^{\ddag\,}$Wu Wen-Tsun Key Laboratory of \vs{-2pt}Mathematics}
\cl{\small University of Science and
Technology of China, Hefei 230026, \vs{2pt}China} \cl{\small E-mail:
ycsu@tongji.edu.cn,\ chgxia@mail.ustc.edu.cn,\
xying@mail.ustc.edu.cn}\vs{5pt}

\par\ni {\small{\bf Abstract.} Intrigued by a well-known theorem of
Mathieu's on Harish-Chandra modules over the Virasoro algebra,
we give  an analogous result for a class of Block type Lie
algebras $\BB$, where the parameter $q$ is a nonzero complex
number. We also classify  quasifinite irreducible highest weight
$\BB$-modules and irreducible $\BB$-modules of the intermediate
series. In particular, we obtain that an irreducible $\BB$-module of
the intermediate series may be a nontrivial extension of a
$\Vir$-module of the intermediate series if
$q$ is half of a negative integer, where $\Vir$ is a subalgebra of $\BB$
isomorphic to the Virasoro algebra. \vskip5pt

\ni{\bf Key words:} Block type Lie algebras; quasifinite modules;
highest weight modules; uniform bounded modules; modules of the
intermediate series.

\ni{\it Mathematics Subject Classification (2000):} 17B10; 17B65;
17B68.}

\vskip10pt\ni {\bf 1. \
Introduction}\setcounter{section}{1}\setcounter{theo}{0}

\ni Since a class of infinite dimensional simple Lie algebras were
introduced by Block [\ref{B}], generalizations of Lie algebras of
this type (usually referred to as {\it Block type Lie algebras})
have been studied by many authors (see, e.g.,
[\ref{DZ},\,\ref{OZ},\,\ref{S3}--\ref{S5},\,\ref{WT}--\ref{X2},\,\ref{ZZ}--\ref{ZM}]).
Even so, the representation theory for Block type Lie algebras is
far from being well developed, except for quasifinite
representations of some particular Block type Lie algebras (see,
e.g., [\ref{S3}--\ref{S5},\,\ref{WT}]). For example, the author of
[\ref{S3},\,\ref{S4}]  studied the representations of the Block type
Lie algebra $\BBB$ with basis
$\{L_{\a,i},c\,|\,\,\a,i\in\Z,i\geq-1\}$ over $\C$ and relations
\begin{equation}\label{B-block}
[L_{\a,i},L_{\b,j}]=(\b(i+1)-\a(j+1))L_{\a+\b,i+j}+
\a\d_{\a+\b,0}\d_{i+j,-2}c, \ \ [c,L_{\a,i}]=0,
\end{equation}
for $\a,\b\in\Z,\,i,j\ge-1$. The author of \cite{WT} studied
representations of the Block type Lie algebra $\BBB(1)$, which can
realized as a special case of Block type Lie algebras considered in
this paper. The author of [\ref{S5}] presented some results on the
classification of quasifinite representations of Lie algebras
related to the Virasoro algebra, including some Block type Lie
algebras.

In this paper, we study systematically representations of Block type
Lie algebras  for a class $\BB$ (not only for a single algebra) with
parameter $q$ being a nonzero complex number, where $\BB$ has  basis
$\{L_{\a,i},c\,|\,\a,i\in\Z,\, i\ge0\}$ over $\C$ and relations
\begin{equation}\label{Bq-block}
[L_{\a,i},L_{\b,j}]=(\b(i+q)-\a(j+q))L_{\a+\b,i+j}
+\d_{\a+\b,0}\d_{i+j,0}\frac{\a^3-\a}{12}c, \ \ [c,L_{\a,i}]=0.
\end{equation}
Note that the Lie algebra $\BBB(0)$ is in fact a half part of the
well-known Virasoro-like algebra, and $\BBB(1)$ is the Block type
Lie algebra studied in \cite{WT}.

The Lie algebra $\BB$ has a natural $\Z$-gradation
$\BB=\oplus_{\a\in\sZ}\BB_\a$ with
\begin{equation}\label{gradeddd}
\BB_\a=\sp\{L_{\a,i}\,|\,\a,i\in\Z,\,i\ge0\}\oplus\d_{\a,0}\C c.
\end{equation}

\begin{defi}\label{defi-module}\rm
\vspace*{-3pt}
\begin{itemize}\parskip-3pt\item[(1)]
A module $V$ over $\BB$ is called\vspace*{-3pt}
\begin{itemize}\parskip-3pt
\item[$\bullet$] {\it $\Z$-graded} if $V=\oplus_{\a\in\Z}V_\a$ and $\BB_\a
  V_\b\subset V_{\a+\b}$ for all $\a,\b$;
\item[$\bullet$] {\it quasifinite} if
  it is $\Z$-graded and ${\rm dim\ssc\,}V_\b <\infty$ for all $\b$;
\item[$\bullet$] {\it uniformly bounded} if it is $\Z$-graded and there is
  $N\!\ge\!0$ with ${\rm dim\ssc\,}V_\b\!\le\!N$ for all $\b$;
\item[$\bullet$] {\it a module of the intermediate series} if it is $\Z$-graded
and ${\rm  dim\ssc\,}V_\b\le1$ for all $\b$;
\item[$\bullet$] a {\it highest}
  (resp.,~{\it lowest}) {\it weight module} if there exists some
  $\L\in\BB_0^*$ (the dual space of $\BB_0$) such that $V=V(\L)$,
  where $V(\L)$ is a module generated by a {\it highest}
  (resp.,~{\it lowest}) {\it weight vector} $v_\L\in V(\L)_0$,
  i.e., $v_\L$ satisfies
$$
hv_\L=\L(h)v_\L\;\;\mbox{for}\;\;\;h\in\BB_0,\;\;\;
\mbox{and}\;\;\;\BB_+v_\L=0 \;\mbox{ (resp.,~$\BB_-v_\L=0$)},
$$
where $\BB_{\pm} =\oplus_{\pm\a>0}\BB_\a$ (cf.~\eqref{Tri---}).
\end{itemize}\item[(2)]A nonzero vector $v$ in a $\Z$-graded module $V$ is called {\it
singular} or {\it primitive} if $\BB_+v\!=\!0$.
\end{itemize}
\end{defi}

When we study representations of a Lie algebra of this kind, as
pointed in [\ref{KL},\,\ref{S1},\,\ref{SX}], we encounter the
difficulty that though it is $\Z$-graded, the graded subspaces are
still infinite dimensional, thus the study of quasifinite modules is
a nontrivial problem. As stated in  \cite{WT}, an important feature
that $\BB$ defined in \eqref{Bq-block} is different from $\BBB$
defined in \eqref{B-block} is that $\BB$ contains the subalgebra
$\Vir$ isomorphic to the well-known Virasoro algebra, where
\begin{eqnarray}\label{denote-Vir}
&\!\!\!\!\!\!\!\!&\Vir=\sp\{L_\a,\kappa\,|\,\a\in\Z\},\ \ \  L_\a:=q^{-1}L_{\a,0},\ \ \kappa:=q^{-2}c,\\
\label{Vir-relat}
&\!\!\!\!\!\!\!\!&
[L_\a,L_\b]=(\b-\a)L_{\a+\b}+\frac{\a^{3}-\a}{12}\delta_{\a+\b,0}\kappa,
\ \ [\kappa,L_\a]=0.
\end{eqnarray}
The authors of [\ref{XYZ}] studied structures of $\ol\BB=\BB/\C c$
with $q$ being a positive integer. We point out here that the
results [\ref{XYZ}, Theorems 2.7, 3.1 and 4.1], including
automorphism groups, derivation algebras and central extensions, of
$\ol\BB$ for $0<q\in\Z$ still hold for
$q\in\C\bs(\frac12{\Z_-}\!\cup\frac13{\Z_-})$, since in case
$q\in\C\bs(\frac12{\Z_-}\!\cup\frac13{\Z_-}\!\cup\Z_+)$, many
coefficients containing $q$ will become invertible. In addition,
$\BB$'s are distinct from each other for different positive rational
number $q$'s,
%
%
%
namely,
$$
\BBB(q_1)\cong\BBB(q_2)\ \ \Longleftrightarrow \ \ q_1=q_2\ \ \
\mbox{for}\ \ \ q_1,q_2\in\Q_+^*.
$$
Furthermore, for any $1<q_1<q_2\in\Z$ with $q_1|q_2$, we find the
following interesting relations:
\begin{equation}\label{Bq-relations}
\BBB(\varepsilon q_2^{-1})\raisebox{-2pt}{$\ ^{\sc\ra}_{\sc\ne}\ $} \BBB(\varepsilon
q_1^{-1})\raisebox{-2pt}{$\ ^{\sc\ra}_{\sc\ne}\
$}\BBB(\varepsilon)\raisebox{-2pt}{$\ ^{\sc\ra}_{\sc\ne}\ $}\BBB(\varepsilon
q_1)\raisebox{-2pt}{$\ ^{\sc\ra}_{\sc\ne}\ $}\BBB(\varepsilon q_2),\mbox{ where }\varepsilon=\pm1.
\end{equation}
More precisely, $\BBB(\varepsilon)$ contains the subalgebra with basis
$\{q_1^{-1}L_{\a,q_1i}\,|\,\a\in\Z, i\in\Z_+\}$ isomorphic to
$\BBB(\varepsilon q_1^{-1})$, and $\BBB(\varepsilon q_2)$ contains the subalgebra
with basis $\{q_2^{-1}L_{\a,q_2 i}\,|\,\a\in\Z, i\in\Z_+\}$
isomorphic to $\BBB(\varepsilon)$.

Moreover, $\BBB(q),\,\BBB(q^{-1})$ with $0<q\in\Z$ are related to
the well-known \mbox{$W$-infinity} Lie algebra $\WW_\infty$ in the
following way: Recall that the $W$-infinity Lie algebra
$\WW_{1+\infty}$ is defined to be the universal central extension of
infinite dimensional Lie algebra of differential operators on the
circle, which has  basis $\{ x^\a D^i, c\,|\,\a\in\Z,i\geq0\}$ with
$D=\frac{d}{dx}$, and relations
$$
[x^\a D^i,x^\b D^j]=x^{\a+\b}((D+\b)^iD^j-D^i(D+\a)^j)+
\d_{\a+\b,0}(-1)^i i! j! \binom{\a+i}{i+j+1}c.
$$
Then the $W$-infinity algebra $\WW_\infty$, the universal central
extension of infinite dimensional Lie algebra of differential
operators on the circle of degree at least one, is simply the
subalgebra of $\WW_{1+\infty}$ spanned by $\{ x^\a
D^i,c\,|\,\a\in\Z,i\geq1\}$. If we define a natural filtration of
$\WW_\infty$ by
$$\begin{array}{ll}
\{0\}=(\WW_\infty)_{[-2]}\subset(\WW_\infty)_{[-1]}\subset\cdots\subset
\WW_\infty,\mbox{ where}\\[7pt]
(\WW_\infty)_{[-1]}=\C c, \mbox{ \ }
(\WW_\infty)_{[n]} = \sp\{x^\a D^i,c\,|\,\a\in\Z,\, 1\le i\le
n+1\}\mbox{ for }n\ge0,
\end{array}$$
then $\BBB(1)$ is simply the associated graded Lie algebra of the
filtered Lie algebra $\WW_\infty$. So roughly speaking, $\BB$
contains (reps., $\BBB(q^{-1})$ lies in)  the associated graded Lie
algebra of $\WW_\infty$ by the fact \eqref{Bq-relations}. As stated
in [\ref{S1}, \ref{S3}, \ref{SX}], the $W$-infinity algebras arise
naturally in various physical theories, such as conformal field
theory, the theory of the quantum Hall effect, etc.; among them the
$\WW_{\infty}$ algebra and $\WW_{1+\infty}$ algebra, of interest to
both mathematicians and physicists, have received intensive studies
in the literature. Due to the importance of the $W$-infinity algebra
$\WW_\infty$, motivated by \eqref{Bq-relations}, it is very natural
to post the following question, which seems to be interesting to us.

\begin{ques}\label{ques1}
Do there exist Lie algebras $\WW(q)$ for $q\in\Z_+$ or
$q^{-1}\in\Z_+$ such that $\WW(1)=\WW_\infty$ and  $
\WW(q_2^{-1})\raisebox{-2pt}{$\ ^{\sc\ra}_{\sc\ne}\ $}
\WW(q_1^{-1})\raisebox{-2pt}{$\ ^{\sc\ra}_{\sc\ne}\
$}\WW(1)\raisebox{-2pt}{$\ ^{\sc\ra}_{\sc\ne}\
$}\WW(q_1)\raisebox{-2pt}{$\ ^{\sc\ra}_{\sc\ne}\ $}\WW(q_2)$ for
$1<q_1<q_2\in\Z$ and $q_1|q_2$.
\end{ques}

Because of the facts stated in the statements before
\eqref{denote-Vir}, one may expect richer results in the
representation theory of $\BB$. Motivated by a well-known result of
Mathieu's in [\ref{M}] (see also \cite{MP,S0}), it is very natural
to consider the classification of quasifinite irreducible
$\BB$-modules. Our first main result is the following theorem (the
analogous results to this theorem for the Virasoro algebra, higher
rank Virasoro algebras, $W$-infinity algebras, and some Block type
Lie algebras  were obtained in
[\ref{LZ},\,\ref{M},\,\ref{S1}--\ref{S4},\,\ref{SX},\,\ref{WT}]).

\begin{theo}\label{MainTheo}
A quasifinite irreducible $\BB$-module is either a highest/lowest
weight module, or a uniformly bounded module.\end{theo}

Our second main result is to give a classification of quasifinite
irreducible highest weight modules. To state the result, we need to
introduce the generating series: For any function $\L\in\BB_0^*$
$($the dual of $\BB_0{\ssc\,})$, we set {\it labels}
$\L_i=\L(L_{0,i})$ for $i\ge0$, and define the following {\it
generating series} with variable $z$,
\begin{equation}\label{Gennnn}
\D_\L(z,q)=2q\mbox{$\sum\limits_{i=0}^\infty$}\frac{z^i}{i!}\L_i+
\mbox{$\sum\limits_{i=0}^\infty$}\frac{z^{i+1}}{i!}\L_{i+1}=\L((2q+zt)t^q
e^{zt}).
\end{equation}
We would like to mention that our generating series $\D_\L(z,1)$
corresponding to the Lie algebra $\BBB(1)$ is different from that
given in [\ref{WT}]. Then our second main result below also recovers
those stated in [\ref{WT}] for quasifinite irreducible
$\BBB(1)$-modules.

\begin{theo}\label{QHWM}
Let $L(\L)$ be an irreducible highest weight module over $\BB$ with
highest weight $\L\in\BB_0^*$. Then  $L(\L)$ is quasifinite if and
only if $\D_\L(z,q)$ is a quasipolynomial.
\end{theo}

Our final main result is to give a classification of irreducible
modules of the intermediate series. To state the result, let us
recall (e.g., \cite{M,MP,S0,S2}) that an indecomposable module of
the intermediate series over $\Vir$ is one of $A_{a,b},\, A_a,\,
B_a,\, a, b\in\C$, or their subquotients, where $A_{a,b},\,
A_a,\,B_a$ all have a basis $\{v_\mu\,|\,\mu\in\Z\}$ with the
trivial action of $c$ and
\begin{eqnarray}
&&A_{a,b}:\ L_{\a,0} v_\mu=q(a+\mu+b\a)v_{\a+\mu},\label{indecomp-case1}\\
&&A_a:\ L_{\a,0} v_\mu=q(\mu+\a)v_{\a+\mu}\ (\mu\neq0),\ \
L_{\a,0} v_0=q\a(a+\a)v_\a,\label{indecomp-case2}\\
&&B_a:\ L_{\a,0} v_\mu=q\mu v_{\a+\mu}\ (\mu\neq-\a),\ \ L_{\a,0}
v_{-\a}=-q\a(a+\a)v_0,\label{indecomp-case3}
\end{eqnarray}
for $\a,\mu\in\Z$ (note that we have the factor ``\,$q$\,'' on the
right-hand sides because of notation \eqref{denote-Vir}). We use
$A'_{0,1}$ to denote the nontrivial subquotients of $A_{0,1}$. Then
a nontrivial irreducible $\Vir$-module of the intermediate series is
isomorphic to either $A'_{0,1}$ or $A_{a,b}$ $(a\notin\Z$ or $b\neq
0,1)$ with the trivial action of $c$ and \vskip-10pt
\begin{eqnarray}
&&A'_{0,1}=\sp\{v_\mu\,|\,\mu\in\Z^*\}:\ L_{\a,0}
v_\mu=q(\mu+\a)v_{\a+\mu},\label{simple-case1}\\
&&A_{a,b}=\sp\{v_\mu\,|\,\mu\in\Z\}:\
L_{\a,0}v_\mu=q(a+\mu+b\a)v_{\a+\mu}.\label{simple-case2}
\end{eqnarray}
Obviously $A'_{0,1}$ or $A_{a,b}$ is also an irreducible
$\BB$-module of the intermediate series (still denoted by $A'_{0,1}$
or $A_{a,b}$) by extending the actions of $L_{\a,i}$ with $i\ge1$
trivially, namely
\begin{equation}\label{tri----}
L_{\a,i}v_\mu=0\mbox{ \ for \ }\a,i\in\Z,\,i\ge1.
\end{equation}
If $q\in\frac{1}{2}\Z_-^*$, for any $s\in\C$, by replacing the
actions \eqref{tri----} by
\begin{equation}\label{simple-extend}
L_{\a,i}v_\mu=
\bigg\{\begin{array}{ll}s v_\mu&\mbox{if \ }(\a,i)=(0,-2q),\\[8pt]
0&\mbox{if \ }i\ge1\mbox{ and }(\a,i)\ne(0,-2q),\end{array}
\end{equation}
we obtain an irreducible $\BB$-module of the intermediate series,
denoted by $A'_{0,1}(s)$ or $A_{a,b}(s)$. Furthermore, if $q=-1$,
for any $s,t\in\C$, the $\Vir$-module $A_{a,b}$ can be defined as a
$\BBB(-1)$-module, denoted by $A_{a,b}(s,t)$, by replacing
\eqref{tri----} by
\begin{equation}\label{simple-extend-q=-1-irr}
L_{\a,i}v_\mu= \left\{\begin{array}{ll}
s v_\mu &\mbox{if \ }(\a,i)=(0,2),\\[8pt]
t v_{\a+\mu} &\mbox{if \ }i=1,\\[8pt]
0 &\mbox{if \ }i\ge2\mbox{ and }(\a,i)\ne(0,2).\end{array}\right.
\end{equation}

Now we can state our final main result below.
\begin{theo}\label{ISM}
Let $V$ be an irreducible $\BB$-module of the intermediate series
such that it is nontrivial as a
$\Vir$-module.\vs{-7pt}\begin{itemize}\parskip-3pt
\item[\rm(1)] If $q\notin\frac12\Z^*_-$, then $V\cong A'_{0,1}$ or $A_{a,b}$ $(a\notin\Z$ or $b\neq 0,1)$.
\item[\rm(2)] If $q\in\frac12\Z^*_-\backslash\{-1\}$, then $V\cong A'_{0,1}(s)$ or $A_{a,b}(s)$ $(a\notin\Z$ or $b\neq 0,1)$,  $s\in\C$.
\item[\rm(3)] If $q=-1$, then $V\cong A'_{0,1}(s)$ or $A_{a,b}(s,t)$ $(a\notin\Z$ or $b\neq 0,1)$,  $s,t\in\C$.
\end{itemize}
\end{theo}


Thus in particular, one sees that an irreducible $\BB$-module of
the intermediate series for $q\in\frac{1}{2}\Z^*_-$, different from
others, can be a nontrivial extension of a $\Vir$-module of the
intermediate series.

Based on Theorem \ref{MainTheo} and results stated in \cite{MP}
(see also Proposition \ref{uniformly}), it is very natural to
conjecture that an irreducible uniformly bounded $\BB$-module is a
module of the intermediate series. Namely,

\begin{conj}
A quasifinite irreducible $\BB$-module is either a highest/lowest
weight module,  or a module of the intermediate series. \end{conj}

Throughout the paper, $q$ is always assumed to be a fixed number in
$\C^*$. We use $\C,$ $\C^*,\Z,\Z^*,\Z_+,\Z_+^*,\Z_-$ and $\Z_-^*$ to
denote respectively the sets of complex, nonzero complex numbers,
integers, nonzero, nonnegative, positive, nonpositive and negative
integers.\vskip15pt

\ni {\bf 2. \ Proof of Theorem
\ref{MainTheo}}\setcounter{section}{2}\setcounter{theo}{0}\setcounter{equation}{0}

\ni We can realize the Lie algebra $\BB$ in the space
$\C[x,x^{-1}]\otimes t^q\C[t]\oplus\C c$ with the bracket
\begin{equation}\label{Bq-block2}{\ssc\!}
[x^\a {\ssc\!}f({\ssc\!}t{\ssc\!}),x^\b {\ssc\!}g({\ssc\!}t{\ssc\!})]\!=\!x^{\a+\b}t^{1-q}(\b f'({\ssc\!}t{\ssc\!})g({\ssc\!}t{\ssc\!})\!-\!\a
f({\ssc\!}t{\ssc\!})g'({\ssc\!}t{\ssc\!})) \!+\!\d_{\a+\b,0}\frac{\a^3\!-\!\a}{12}{\rm
Res}_t(t^{-2q-1}\!{\ssc\!}f({\ssc\!}t{\ssc\!})g({\ssc\!}t{\ssc\!}))c,\!\!
\end{equation}
for $\a,\b\in\Z$ and $f(t),g(t)\in t^q\C[t]$, where the prime stands
for the derivative ${d\over dt}$, and ${\rm Res}_tf(t)$ stands for
the {\it residue} of $f(t)$, namely the coefficient of $t^{-1}$ in
$f(t)$. We always denote
\begin{equation}\label{L-a-i}
L_{\a,i}=x^\a t^{q+i}.
\end{equation}
Thus (\ref{Bq-block2}) is equivalent to
(\ref{Bq-block}). Using the gradation \eqref{gradeddd}, we introduce the following notations for $\b,\g\in\Z$,
$$
\BB_{[\b,\g]}=\mbox{$\sum\limits_{\b\le\a\le\g}$}\BB_\a,
$$
and similarly for $\BB_{[\b,+\infty)}$, $\BB_{[\b,\g)}$ and so on.
Putting $\BB_{\pm} =\oplus_{\pm\a>0}\BB_\a$, we have the following
triangular decomposition:
\begin{equation}\label{Tri---}
\BB=\BB_-\oplus\BB_0\oplus\BB_+.
\end{equation}
Note that $\BB_0=t^q\C[t]\oplus\C c$ is an infinite dimensional
commutative subalgebra of $\BB$ (but not a Cartan subalgebra).

Now suppose $V=\oplus_{\mu\in\Z}V_\mu$ is a quasifinite
$\BB$-module. Taking $\mu_0\in\Z^*$, since $c|_{V_{\mu_0}}$ (the
action of $c$ on $V_{\mu_0}$) and $t^{q+i}|_{V_{\mu_0}} $ for
$i\in\Z_+$ are linear transformations of the finite dimensional
subspace $V_{\mu_0}$, there exists big enough fixed integer $p_0$
 such that the operators
$c|_{V_{\mu_0}}, t^q|_{V_{\mu_0}}, \ldots, t^{q+p-1}|_{V_{\mu_0}}$
are linear dependent for all $p\ge p_0$. Therefore, for any $p\ge
p_0$, there exists $f_p(t)\in\BB_0$ of degree $q+p-1$ and $m_p\in\C$
such that
\begin{equation}\label{monomial}
(f_p(t)+m_p c)v=0\ \ \mbox{for}\ \ v\in V_{\mu_0}.
\end{equation}
Define the Lie subalgebra $\LL$ of $\BB$ as follows
\begin{equation*}
\LL=\left\{
\begin{aligned}
&\left\langle x^{-\mu_0}t^q,x^{-\mu_0}t^{q+1},
x^{-\mu_0}t^{q+2},x^{-\mu_0+1}t^q,
f_p(t)+m_p c\,|\,p\ge p_0\right\rangle\ \ \mbox{if}\ \ \mu_0\le-1,\\
&\left\langle x^{-\mu_0}t^q,x^{-\mu_0}t^{q+1},
x^{-\mu_0}t^{q+2},x^{-\mu_0-1}t^q, f_p(t)+m_p c\,|\,p\ge
p_0\right\rangle\ \ \mbox{if}\ \ \mu_0\ge 1,
\end{aligned}\right.
\end{equation*}
where the angle bracket $\langle\ ,\rangle$ stands for ``the Lie
subalgebra generated by''.

\begin{lemm}\label{Lemm-2-1} For any $s\ge1$, and fixed $\mu_0\in\Z^*$,
we hav\vspace*{-3pt}e \baselineskip3pt\lineskip7pt\parskip-3pt
\begin{itemize}\parskip-3pt
  \item[{\rm(1)}] if $\mu_0\le-1$, then there exists $\a_s\in\Z_+^*$
  such that $x^\a t^{q+s-1}\!\in\!\LL$ for all $\a\ge\a_s$;
  \item[{\rm(2)}] if $\mu_0\ge 1$, then there exists $\a_s\in\Z_-^*$
  such that $x^\a t^{q+s-1}\in\LL$ for all $\a\le\a_s$.
\end{itemize}
\end{lemm}
\ni{\it Proof.~}~We only prove part (1) by induction on $s$ (part
(2) can be proved similarly). In case $s=1$, Remark \ref{lower
bound} below shows that, for any integer $\a\ge (1-\mu_0)^2$, there
exist two positive integers $k_1, k_2$ such that
\begin{equation}\label{a-k1-k2}\a=k_1(1-\mu_0)-k_2\mu_0.\end{equation}
 Letting $z_1=x^{-\mu_0+1}t^q$,
$z_2=x^{-\mu_0}t^q$, using \eqref{Bq-block2} and by induction on
$k_1,k_2$, we obtain
\begin{equation}\label{s====1}
{\rm ad}_{z_2}^{k_2-1}{\rm
ad}_{z_1}^{k_1}(z_2)=q^{k_1+k_2-1}\mbox{$\prod\limits_{i=1}^{k_1}
(-(i-1)\mu_0+i-2)\prod\limits_{j=1}^{k_2-1}$}(-(k_1+j-1)\mu_0+k_1)x^\a
t^q.
\end{equation}
Note that the coefficient of $x^\a t^q$ on the right-hand side of
\eqref{s====1} is nonzero. Hence, $x^\a t^q\in\LL$. Now suppose
$s>1$, and inductively assume that there exists an integer
$\a_{s-1}$ such that $x^\a t^{q+s-2}\in\LL$ for $\a\ge\a_{s-1}$. We
denote $r_{\a,q}=\a$ if $s=3,\,q=-1$, or else $r_{\a,q}
=\mu_0(2q+s-1)+\a(q+1)$. We can always choose big enough $\a_s'$
such that $r_{\a,q}\ne0$ whenever $\a\ge\a_s'$. Now take
$\a_s=\max\left\{\a_{s-1}-\mu_0,\a_s'\right\}$, then for all
$\a\ge\a_s$, we have
$$
x^\a t^{q+s-1}=\left\{\begin{array}{lll}
\dis-\frac{1}{r_{\a,q}}[x^{\a+\mu_0} t^{q+s-3},x^{-\mu_0}t^{q+2}]
&\mbox{if \ }s=3,\,q=-1,\\[12pt]
\dis-\frac{1}{r_{\a,q}}[x^{\a+\mu_0} t^{q+m-2},x^{-\mu_0}
t^{q+1}]&\mbox{else},
\end{array}\right.$$ which shows $x^\a t^{q+s-1}\in\LL$. Part
(1) is proved. \qed

\begin{rema}\label{lower bound}\rm
The lower bound $(1-\mu_0)^2$ of $\a$ for case $s=1$ in the above
lemma is more precise than that for a quasifinite $\BBB(1)$-module
given in [\ref{WT}], 
which can be deduced as follows: For any $\a\ge(1-\mu_0)^2$, we
denote $k_0=[\frac{\a}{1-\mu_0}]$ (the integral part of
$\frac{\a}{1-\mu_0}$). Then $k_0\ge1-\mu_0$, and two integers
$$\begin{array}{lll}
k_1:=\a+(k_0+1)\mu_0\ge
k_0(1-\mu_0)+(k_0+1)\mu_0=k_0+\mu_0\ge1,\\[4pt]
k_2:=(k_0+1)(1-\mu_0)-\a>0,
\end{array}$$
satisfy \eqref{a-k1-k2}.
\end{rema}

\begin{lemm}\label{Lemm-2-2} Let $V=\oplus_{\mu\in\Z}V_\mu$ be a
quasifinite irreducible $\BB$-module.
\baselineskip3pt\lineskip7pt\parskip-3pt
\begin{itemize}\parskip-3pt
  \item[{\rm(1)}] If $\mu_0\le-2$, and there exists
  $0\ne v_0\in V_{\mu_0}$ satisfying $\BB_{[\a,+\infty)}v_0=0$
  for some $\a>0$, then $V$ has a highest weight vector.
  \item[{\rm(2)}] If $\mu_0\ge 2$, and there exists
  $0\ne v_0\in V_{\mu_0}$ satisfying $\BB_{(-\infty,\a]}v_0=0$
  for some $\a<0$, then $V$ has a lowest weight vector.
\end{itemize}
\end{lemm}

One can prove this lemma in a similar way as in \cite{S1,S3,WT}, and
the details are omitted.
\vskip5pt

\ni{\it Proof of Theorem \ref{MainTheo}.}\ \ Assume that
$V=\oplus_{\mu\in\Z}V_\mu$ is a quasifinite irreducible $\BB$-module
without highest and lowest weight vectors. We should prove that
\begin{equation}
\dim V_\mu\le\left\{
\begin{aligned}
&3{\ssc\,}\dim V_0+\dim V_1 &&\mbox{if}\ \ \mu\le-2,\\
&3{\ssc\,}\dim V_0+\dim V_{-1} &&\mbox{if}\ \ \mu\ge 2.
\end{aligned}\right.
\end{equation}
For fixed $\mu_0\le-2$, we claim that the following linear map is
injective:
$$
\T_{\mu_0}^-=\Big(x^{-\mu_0}t^q\oplus x^{-\mu_0}t^{q+1}\oplus
x^{-\mu_0}t^{q+2}\oplus x^{-\mu_0+1}t^q\Big)\Big|_{V_{\mu_0}}: V_{\mu_0}\rightarrow
V_0\oplus V_0\oplus V_0\oplus V_1.
$$
Otherwise there exists $0\ne v_0\in V_{\mu_0}$ such that
$\T_{\mu_0}^-(v_0)=0$, which implies that $x^{-\mu_0}t^q$,
$x^{-\mu_0}t^{q+1}$, $x^{-\mu_0}t^{q+2}$ and $x^{-\mu_0+1}t^q$ take
$v_0$ to zero. On the other hand,  $(f_p(t)+m_p c)v_0=0$ for
$p\ge p_0$ by \eqref{monomial}. Hence, by definition,
\begin{equation}\label{subalg-kill}
\LL v_0=0.
\end{equation}
Applying Lemma \ref{Lemm-2-1}(1), for any $1\le p< p_0$, there
exists some positive integer $\a_p$ such that $x^\a t^{q+p-1}\in\LL$
for $\a\ge\a_p$. Denote $\G=\max\{\a_1,\a_2,\ldots,\a_{p_0-1}\}$. Then $x^\a t^{q+p-1}\in\LL$ for $1\le p< p_0$, $\a\ge\G$.
Furthermore, for $p\ge p_0$, $\a\ge\G$, we have
\begin{equation}\label{case-p>p0}
x^\a tf'_p(t)=-\frac{1}{\a}[x^\a t^q,f_p(t)+m_p c]\in\LL.
\end{equation}
Taking $p=p_0+i$ in \eqref{case-p>p0}, noting that $\deg
f_p(t)=q+p-1$, we have $x^\a t^{q+p_0+i-1}\in\LL$ for $\a\ge\G$,
$i\in\Z_+$. Therefore $x^\a t^{q+p-1}\in\LL$ for $p\ge1$, $\a\ge\G$,
namely,
\begin{equation}\label{subalg-contain}
\BB_{[\G,+\infty)}\subseteq\LL.
\end{equation}
By \eqref{subalg-kill} and \eqref{subalg-contain},
$\BB_{[\G,+\infty)}v_0=0$. Then Lemma \ref{Lemm-2-2}(1) shows $V$
has a highest weight vector, which contradicts our assumption. Thus
the map $\T_{\mu_0}^-$ is injective, which implies $\dim V_\mu\le
3{\ssc\,}\dim V_0+\dim V_1$ if $\mu\le-2$.

Similarly, one can derive $\dim V_\mu\le 3{\ssc\,}\dim V_0+\dim
V_{-1}$ if $\mu\ge 2$ by Lemma \ref{Lemm-2-1}(2) and Lemma
\ref{Lemm-2-2}(2). Denote $N=\max\{3{\ssc\,}\dim V_0+\dim
V_1,3{\ssc\,}\dim V_0+\dim V_{-1}\}$. Then $\dim V_\mu\le N$ for
$\mu\in\Z$, namely $V$ is a uniformly bounded $\BB$-module. This
completes the proof.\qed\vskip15pt

\ni {\bf 3. \ Quasifinite highest weight
modules}\setcounter{section}{3}\setcounter{theo}{0}\setcounter{equation}{0}

\ni In this section,  we start with general settings on parabolic
subalgebras of $\Z$-graded Lie algebra. Following
\cite{KL,KR,S3,S4}, after giving some descriptions of parabolic
subalgebras of $\BB$, we use the results to characterize the
irreducible quasifinite highest weight $\BB$-modules by generating
series.

\begin{defi}\label{defi-parabolic}\rm
Let $\LLL=\oplus_{\a\in\Z}\LLL_\a$ be a $\Z$-graded Lie algebra.
\baselineskip3pt\lineskip7pt\parskip-3pt
\begin{itemize}\parskip-3pt
\item[(1)]
A subalgebra $\PPP$ of $\LLL$ is called {\it parabolic} if it
contains $\LLL_0\oplus\LLL_+$ as a proper subalgebras, namely, $\PPP=\oplus_{\a\in\Z}\PPP_\a$ with $\PPP_\a=\LLL_\a$ for
$\a\ge0$, and $\PPP_\a\ne\{0\}$ for some $\a<0$.
\item[(2)]
Given $0\ne a\in\LLL_{-1}$, we define a  parabolic subalgebra
$\PPP(a)=\oplus_{\a\in\Z}\PPP(a)_\a$ of $\LLL$ 
as follows:
\begin{equation}\label{defi-mini-para}
\PPP(a)_\a=\left\{
\begin{aligned}
&\LLL_\a &&\mbox{if}\ \ \a\ge0,\\
&\sp\{[\ldots,[\LLL_0,[\LLL_0,a]]\cdots]\} &&\mbox{if}\ \ \a=-1,\\
&[\PPP(a)_{-1}, \PPP(a)_{\a+1}] &&\mbox{if}\ \ \a\le-2.
\end{aligned}\right.
\end{equation}
\item[(3)]
A parabolic subalgebra $\PPP$ is called {\it nondegenerate} if
$\PPP_\a$ has finite codimension in $\LLL_\a$ for all $\a<0$.\item[(4)]A
nonzero element $a\in\LLL_{-1}$ is called {\it nondegenerate} if
$\PPP(a)$ is nondegenerate.
\end{itemize}\end{defi}

A {\it Verma module} over $\BB$ is defined as the induced module
$$
M(\L)=U(\BB)\otimes_{U(\BB_0\oplus\BB_+)}\C v_\L\ \ \mbox{for}\ \
\L\in\BB_0^*,
$$
where $\C v_\L$ is the one-dimensional $\BB_0\oplus\BB_+$-module
given by $(h+n)(v_\L)=\L(h)v_\L$ for $h\in\BB_0, n\in\BB_+$ (cf.~\eqref{Tri---}). Here
and further $U(\LLL)$ stands for the universal enveloping algebra of
a Lie algebra $\LLL$. Then any highest weight module $V(\L)$ is a
quotient module of $M(\L)$ and the irreducible highest weight module
$L(\L)$ is the quotient of $M(\L)$ by the maximal proper $\Z$-graded
submodule.

Define a parabolic subalgebra
$\PPP(q,a)=\oplus_{\a\in\Z}\PPP(q,a)_\a$ of $\BB$ as in
\eqref{defi-mini-para}, where $0\ne a\in\BB_{-1}$. By [\ref{KL},
Lemma 2.2], $\PPP(q,a)$ is the minimal parabolic subalgebra
containing $a$ and
$$
\BBB(q,a)_0:=\BB_0\cap[\PPP(q,a),\PPP(q,a)]=[a,\BB_1].
$$
Write $a=x^{-1}f(t)$, then we have
$[a,xg(t)]=[x^{-1}f(t),xg(t)]=(f'(t)g(t)+f(t)g'(t))t^{1-q}$ for
$g(t)\in t^q\C[t]$, which implies
\begin{equation}\label{bqa-span}
\BBB(q,a)_0=\sp\{(f(t)g(t))'t^{1-q}\,|\,g(t)\in t^q\C[t]\}.
\end{equation}
Let $\L\in\BB_0^*$ be such that $\L|_{\BB_0\cap[\PP,\PP]}=0$. Then
the $\BB_0\oplus\BB_+$-module $\C v_\L$ can be extended to  be a
$\PP$-module by letting $\PP_\a$ take $v_\L$ to zero for $\a<0$. We
construct the following highest weight $\BB$-module
$$
M(\PP,\L)=U(\BB)\otimes_{U(\PP)}\C v_\L,
$$
which is called the {\it generalized Verma module}.

\begin{lemm}\label{lemm-3-1}
Let $\PP=\oplus_{\a\in\Z}\PP_\a$ be a parabolic subalgebra of $\BB$.
\baselineskip3pt\lineskip7pt\parskip-3pt
\begin{itemize}\parskip-3pt
\item[{\rm(1)}] There exists an nonzero element $0\ne a\in\BB_{-1}$ such that
  $\PPP(q,a)\subseteq\PP$.
\item[{\rm(2)}] For any $\a<0$, the subspace $\PP_\a$ is nontrivial, and has
  finite codimension in $\BB_\a$.
\item[{\rm(3)}] $\PP$ is nondegenerate, and any nonzero element
  $0\ne a\in\BB_{-1}$ is nondegenerate.
\end{itemize}
\end{lemm}

\ni{\it Proof.} \ \ (1) By definition, there exists at least one
$\a<0$ such that $\PP_\a\ne\{0\}$. We claim that
$\PP_{\a+1}\ne\{0\}$ if $\a\le-2$. Otherwise $[\PP_\a,\BB_1]=0$.
Since $\a<0$, we can easily choose some positive integer $j_0$ such
that $k_\a=(q+i)-\a(q+j_0)\ne0$ for $i\in\Z_+$. Taking any $0\ne
b=\sum_{i\in I}b_i x^\a t^{q+i}\in\PP_\a$, where $I$ is a finite
subset of $\Z_+$ and $b_i\in\C$, we have
\begin{equation}\label{fact1}
0=[b,xt^{q+j_0}]=\Big[\mbox{$\sum\limits_{i\in I}$}b_i x^\a
t^{q+i},xt^{q+j_0}\Big]=\mbox{$\sum\limits_{i\in I}$}b_i k_\a
x^{\a+1}t^{q+i+j_0},
\end{equation}
which implies $b_i=0$ for $i\in I$,  i.e., $b=0$, a contradiction.
This proves the claim. Therefore $\PP_{-1}\ne\{0\}$ by induction.
Taking any nonzero element $a\in\PP_{-1}$, we have
$\PPP(q,a)\subseteq\PP$ by the minimality of $\PPP(q,a)$.

(2) We shall use induction on $\a<0$ to show $\PP_\a\ne\{0\}$. The
case $\a=-1$ is proved in (1). Now suppose $\PP_\a\ne\{0\}$ for some
$\a\le-1$. For $0\ne x^\a f(t)\in\PP_\a$, we have
\begin{equation}\label{fact2}
x^\a f(t)t^i=\frac{1}{\a(q+i)}[t^{q+i},x^\a f(t)]\in\PP_\a\ \
\mbox{if}\ \ i\ne-q.
\end{equation}
Let $z_1:=x^\a f_1(t)\in\PP_\a$ and $z_2:=x^{-1} f_2(t)\in\PP_{-1}$
be any nonzero elements. Then $z_3:=[z_1,z_2]\in \PP_{\a-1}$. If
$q\ne-1$, then \eqref{fact2} with $i=1$ implies
$z_3t\in\PP_{\a-1},\,z_2t\in\PP_{-1}$, and so $$\a
x^{\a-1}f_1(t)f_2(t)t^{1-q}=z_3t-[z_1,z_2t]\in\PP_{\a-1},$$ which is
clearly a nonzero element. If $q=-1$, then \eqref{fact2} with $i=2$
implies $z_3t^2\in\PP_{\a-1},\,z_2t^2\in\PP_{-1}$, and so
$$
2\a x^{\a-1}f_1(t)f_2(t)t^3=z_3t^2-[z_1,z_2t^2]\in\PP_{\a-1},
$$
which is clearly a nonzero element. Thus by induction,
$\PP_\a\ne\{0\}$ for all $\a<0$. This together with  \eqref{fact2}
immediately implies that $\PP_\a$ has finite codimension in
$\BB_\a$.

(3) By definition, $\PP$ is nondegenerate by (2). In particular,
$\PPP(q,a)$ is nondegenerate for any nonzero element $0\ne
a\in\BB_{-1}$, namely $a$ is nondegenerate by definition.
\qed\vskip10pt

Using \eqref{fact1}, Lemma \ref{lemm-3-1} and [\ref{KL}, Theorem
2.5], we have the following lemma.

\begin{lemm}\label{lemm-3-2}
The following  conditions on $\L\in\BB_0^*$ are equivalent:
\baselineskip3pt\lineskip7pt\parskip-3pt
\begin{itemize}\parskip-3pt
\item[{\rm(1)}] $L(\L)$ is quasifinite;
\item[{\rm(2)}] there exists an element $0\ne a\in\BB_{-1}$ such that
$\L(\BBB(q,a)_0)=0$;
\item[{\rm(3)}] $M(\L)$ contains a singular vector
$a\cdot v_\L\in M(\L)_{-1}$ $($cf.~Definition $\ref{defi-module}(2))$, where $0\ne a\in\BB_{-1}$;
\item[{\rm(4)}] there exists an element $0\ne a\in\BB_{-1}$ such that
$L(\L)$ is an irreducible quotient of the generalized Verma module
$M(\PPP(q,a),\L)$.
\end{itemize}
\end{lemm}

Assume that $L(\L)$ is a quasifinite irreducible highest weight
module over $\BB$. By Lemma \ref{lemm-3-2}, there exists some monic
polynomial $f(t)\in t^q\C[t]$ such that $(x^{-1}f(t))v_{\L}=0$. We
shall call such monic polynomial of minimal degree, uniquely
determined by the highest weight $\L$, the {\it characteristic
polynomial} of $L(\L)$.

Recall that a {\it quasipolynomial} is a linear combination of
functions of the form $p(z)e^{az}$, where $p(z)$ is a polynomial and
$a\in\C$. A well-known fact
[\ref{KL},\,\ref{KR},\,\ref{S3}--\ref{S5}] stated that a formal
power series is a quasipolynomial if and only if it satisfies a
nontrivial linear differential equation with constant coefficients.\vs{5pt}

\ni{\it Proof of Theorem \ref{QHWM}.}~~Clearly, $f(t)e^{zt}=f(\frac{\p}{\p z})e^{zt}$
for $f(t)\in\C[t]$, here and further we use notation $
e^{zt}=\mbox{$\sum_{i=0}^\infty$}\frac{z^i}{i!}t^i $ as a generating
series of $\C[t]$. For any $f(t)\in t^q\C[t]$, we denote $\wt
f(t):=t^{-q}f(t)\in\C[t]$, then $f(t)e^{zt}=\wt f(\frac{\p}{\p
z})(t^qe^{zt})$. Recalling that the prime stands for the partial
derivative $\frac{\p}{\p t}$, we have
\begin{eqnarray}\label{quasipolynomial}
\L((f(t)t^q e^{zt})'t^{1-q})&=&\L\Big(\big(\wt f(\mbox{$\frac{\p}{\p
z}$})(t^{2q}e^{zt})\big)'t^{1-q}\Big)\nonumber\\
&=&\Big(\big(\wt f(\mbox{$\frac{\p}{\p z}$})(2q
t^{2q-1}e^{zt}+z t^{2q}e^{zt})\big)t^{1-q}\Big)\nonumber\\
&=&\wt f(\mbox{$\frac{\p}{\p z}$})\L((2q+zt)t^q e^{zt})=\wt
f(\mbox{$\frac{d}{d z}$})\D_\L(z,q).
\end{eqnarray}
If $L(\L)$ is quasifinite, then by \eqref{bqa-span} and Lemma
\ref{lemm-3-2}(2) there exists a polynomial $0\ne f(t)\in t^q\C[t]$
such that $\L((f(t)g(t))'t^{1-q})=0$ for all $g(t)\in t^q\C[t]$.
Taking $g(t)=t^q e^{zt}$, by \eqref{quasipolynomial}, we have $\wt
f(\frac{d}{d z})\D_\L(z,q)=0$, which implies that $\D_\L(z,q)$ is a
quasipolynomial.

Conversely, if $\D_\L(z,q)$ is a quasipolynomial, then there exists
a polynomial $0\ne h(t)\in\C[t]$ such that $h(\frac{d}{d
z})\D_\L(z,q)=0$. Denote $f(t)=t^q h(t)\in t^q\C[t]$, then $\wt
f(\frac{d}{d z})\D_\L(z,q)=0$. By \eqref{quasipolynomial}, we have
$$
0=\L((f(t)t^q e^{zt})'t^{1-q})
=\L\left(\Big(f(t)\mbox{$\sum\limits_{i=0}^\infty$}\frac{z^i}{i!}t^{q+i}\Big)'t^{1-q}\right)
=\mbox{$\sum\limits_{i=0}^\infty$}\frac{z^i}{i!}\L\big((f(t)t^{q+i})'t^{1-q}\big),
$$
which implies $\L((f(t)t^{q+i})'t^{1-q})=0$ for $i\in\Z_+$. Hence
$\L((f(t)g(t))'t^{1-q})=0$ for $g(t)\in t^q\C[t]$ and thus $L(\L)$
is quasifinite by \eqref{bqa-span} and Lemma \ref{lemm-3-2}(2).
\qed\vskip15pt

\ni {\bf 4. \ Intermediate series
modules}\setcounter{section}{4}\setcounter{theo}{0}\setcounter{equation}{0}

\ni Suppose $V=\oplus_{\mu\in\Z}V_\mu$ is an irreducible uniformly
bounded $\BB$-module which is nontrivial as a $\Vir$-module. For any $a\in\C$, we let
$$
V[a]=\mbox{$\bigoplus\limits_{\mu\in\Z}$}V_\mu[a],\ \ \mbox{where}\
\ V_\mu[a]=\{v\in V_\mu\,|\,L_{0,0}v=q(\mu+a)v\}.
$$
By \eqref{Bq-block}, one can check that $V[a]$ is a $\BB$-submodule,
which is a direct summand of $V$. Thus $V=V[a]$ for some fixed
$a\in\C$, namely,
\begin{equation}\label{V-module}
V=\mbox{$\bigoplus\limits_{\mu\in\Z}$}V_\mu,\ \ \mbox{where}\ \
V_\mu=V_\mu[a].
\end{equation}
Note that, regarding as a $\Vir$-module, $V$ is also uniformly
bounded. Therefore, by the results of \cite{MP,S0, S2}, we have the
following proposition.

\begin{prop}\label{uniformly}If $V$ is an irreducible
uniformly bounded $\BB$-module as in $(\ref{V-module})$, then there
exists a non-negative integer $N$ such that ${\rm
dim\ssc\,}V_\mu[a]=N$ for all $\mu\in\Z$ with $\mu+a\neq0$.
\end{prop}

The following result is well-known
(cf.~\eqref{indecomp-case1}--\eqref{indecomp-case3}).

\begin{lemm}\label{ISM-Vir-redu}
Let $V=\oplus_{\mu\in\Z}V_\mu[a]$ be a reducible $\Vir$-module of
the intermediate series, then $V$ is isomorphic to one of $A_a$,
$B_a$ or $A'_{0,1}\oplus \C v_0$ as a direct sum of $\Vir$-modules.
\end{lemm}

The following
 lemma seems to be crucial in obtaining Theorem \ref{ISM}.
\begin{lemm}\label{Th--1}
An irreducible $\BB$-module of the intermediate series $V$ remains
to be irreducible when regarded as a $\Vir$-module. In particular,
$V$ remains irreducible as a $\BBB(\frac{q}{k})$-module for any
$k\in\Z_+^*$.
\end{lemm}
\ni{\it Proof.~}~The second statement follows from the first since
$\Vir\subset\BBB(\frac{q}{k})$. We prove the first statement in two cases.
\vskip4pt

{\it Case 1}: $q\ne-1$. \vskip4pt

If the statement is not true, then there exists
a proper irreducible $\Vir$-submodule $M$. First suppose $M=M_0$ is
trivial. Then $L_{\a,0}M_0=0,\,L_{0,i}M_0\subset M_0$. Since $\BB$
is generated by $\{L_{\a,0},L_{0,i}\,|\,\a,i\in\Z,\, i\ge1\}$, we
see $M$ is a proper $\BB$-submodule, a contradiction with the
irreducibility of $V$.

Now suppose $M$ is nontrivial, which has to have the form
\eqref{simple-case1}. Thus $M_\mu:=M\cap V_\mu=V_\mu$ if $\mu\ne0$
and $M_0=0$. Then for any $\a,\mu\in\Z$ and $i\in\Z_+$,
\begin{equation}\label{MA1}
L_{0,i}M_0=0,\ \ L_{\a,0}M_{-\a}=0\mbox{ \ and \ }
L_{\a,i}M_\mu\subset V_{\a+\mu}=M_{\a+\mu}\mbox{ \ if \ }\mu\ne-\a.
\end{equation}
Furthermore, for $\a\ne0$, we have
\begin{equation}\label{MA1-2}
\a(q+i)L_{\a,i}M_{-\a}=[L_{\a,0},L_{0,i}]M_{-\a}
=L_{\a,0}L_{0,i}M_{-\a}-L_{0,i}L_{\a,0}M_{-\a}=0,
\end{equation}where the last equality follows from
\eqref{MA1}. Since $q\ne-1$, \eqref{MA1} and \eqref{MA1-2} in
particular imply $L_{\a,i}M\subset M$ for all $\a,i$ with $i\le1$.
Since $\BB$ is generated by $\{L_{\a,i}\,|\,\a,i\in\Z,\,0\le
i\le1\}$, we see $M$ is a nontrivial proper $\BB$-submodule, a
contradiction with the irreducibility of $V$. So, $V$ is an
irreducible $\Vir$-module of the intermediate series.
\vskip4pt

{\it Case 2}: $q=-1$. \vskip4pt

 Suppose $V$ becomes reducible
when regarded as a $\Vir$-module, which is isomorphic to $A_a$,
$B_a$ or $A'_{0,1}\oplus\C v_0$ by Lemma \ref{ISM-Vir-redu}. We should
show the following claim, which
leads to a contradiction.
 \begin{clai}
 $\C v_0$ is a submodule or a quotient module of $V$
 \end{clai}

 Since $\BBB(-1)$ can be generated by
$\{L_{1,1},L_{0,2},L_{\a,0}\,|\,\a\in\Z\}$, it suffices to determine
the actions of $L_{1,1}$ and $L_{0,2}$. Suppose
$L_{\a,1}v_\mu=e_{\a,\mu} v_{\a+\mu}$, $L_{\a,2}v_\mu=f_{\a,\mu}
v_{\a+\mu}$, and write $e_{1,\mu}=e_\mu$,
$f_{0,\mu}=f_\mu$ for short.
\vskip4pt

{\it Subcase 2.1}: $V\cong A_a$ as $\Vir$-modules. \vskip4pt

We prove that $e_\mu=0$ for
$\mu\ne0$, and $f_\mu$ is a constant. Then, together with
\eqref{indecomp-case2}, it is easy to check that $\C v_0$ is a
quotient module of $V$ (thus Claim 1 holds). Applying $L_{\a,1}=[L_{1,1},L_{\a-1,0}]$ to
$v_\mu$, we obtain
\begin{eqnarray*}
e_{\a,\mu}= \left\{\begin{array}{ll}
(\mu+\a)e_\mu-(\mu+\a-1)e_{\mu+\a-1} &\mbox{if \ }\mu\ne0,-1,\\[8pt]
-(\a-1)(a+\a-1)e_{\a-1}+\a e_0 &\mbox{if \ }\mu=0,\\[8pt]
-(\a-2)e_{\a-2}+(\a-1)(a+\a-1)e_{-1} &\mbox{if \
}\mu=-1.\end{array}\right.
\end{eqnarray*}
Applying $[L_{-1,1},L_{1,0}]=-L_{0,1}$ to $v_0$ gives $e_{-1}=0$.
Using this and applying $[L_{\mu,1},L_{1,0}]=\mu L_{\mu+1,1}$ to
$v_{-1}$ with $\mu\ne1$, we obtain
$(\mu-1)e_{\mu-1}=(\mu-2)e_{\mu-2}$, which implies that $e_\mu=0$
for $\mu\ne0$. Similarly, applying
$L_{\a,2}=\frac{1}{\a}[L_{0,2},L_{\a,0}]$ with $\a\ne0$ to $v_\mu$,
we obtain
\begin{eqnarray*}
f_{\a,\mu}= \left\{\begin{array}{ll}
\frac{1}{\a}(\mu+\a)(f_{\mu}-f_{\a+\mu}) &\mbox{if \ }\mu\ne0,\\[8pt]
(a+\a)(f_{0}-f_\a) &\mbox{if \ }\mu=0.\end{array}\right.
\end{eqnarray*}
Applying $[L_{\mu,2},L_{1,0}]=(\mu+1) L_{\mu+1,2}$ to $v_0$ with
$\mu\ne0,-1$, and $v_{-1}$ with $\mu\ne0,\pm1$ respectively, we
obtain
$$\begin{array}{ll}
\Eqa_1(\mu):=\mu(a+\mu)(f_0-f_{\mu})-(a+1)(f_1-f_{\mu+1})-\mu(a+\mu+1)(f_0-f_{\mu+1})=0,\\[6pt]
\Eqa_2(\mu):=(\mu-1)(f_{-1}-f_{\mu-1})-\mu(f_{-1}-f_\mu)=0.\end{array}$$
Solving following linear equations on
$f_0,f_{\pm1},f_{\pm2},f_{\pm3}$:
\begin{eqnarray*}
\Bigg\{\begin{array}{lllll}
\Eqa_1(2)=0,&\Eqa_1(-2)=0,&\Eqa_1(-3)=0,&\Eqa_1(1)=0,\\[8pt]
\Eqa_2(2)=0,&\Eqa_2(-2)=0,&\Eqa_2(-3)=0,
\end{array}
\end{eqnarray*}
we obtain $f_0=f_{\pm1}=f_{\pm2}=f_{\pm3}$. Now, rewriting
$\Eqa_2(\mu)=0$ as $\mu f_\mu-(\mu-1)f_{\mu-1}=f_{-1}$, we see
that $f_\mu$ is a constant. Thus Claim  1 holds.
\vskip4pt

{\it Subcase 2.2}: $V\cong B_a$ as $\Vir$-modules. \vskip4pt

 We prove that $e_\mu=0$ for
$\mu\in\Z$, and $f_\mu$ is a constant. Then, together with
\eqref{indecomp-case3}, it is easy to check that $\C v_0$ is a
submodule of $V$ (thus Claim 1 holds). Applying $L_{\a,1}=[L_{1,1},L_{\a-1,0}]$ to
$v_\mu$, we obtain
\begin{eqnarray*}
e_{\a,\mu}= \left\{\begin{array}{ll}
(\mu+1)e_\mu-\mu e_{\a+\mu-1} &\mbox{if \ }\mu\ne-\a,-\a+1,\\[8pt]
\a e_{-1}-(\a-1)(a+\a-1)e_{-\a} &\mbox{if \ }\mu=-\a,\\[8pt]
(\a-1)(a+\a-1)e_0-(\a-2)e_{-\a+1} &\mbox{if \
}\mu=-\a+1.\end{array}\right.
\end{eqnarray*}
Applying $[L_{\mu,1},L_{1,0}]=\mu L_{\mu+1,1}$ to $v_{-\mu}$ with
$\mu\ne1$, and $v_{-\mu-1}$ with $\mu\ne0$ respectively, we obtain
$$\begin{array}{rl}\Eqb_1(\mu):=\!\!\!&\mu((\mu+1)e_{-2}-\mu
e_{-\mu-1})-\mu(a+\mu)e_0+(\mu-1)e_{-\mu}=0,\\[8pt]
\Eqb_2(\mu):=\!\!\!&(\mu+1)(\mu e_{-1}-(\mu-1)(a+\mu-1)e_{-\mu})-
(a+1)((\mu+1)e_{-2}-\mu
e_{-\mu-1})\\[5pt]&-\mu(\mu(a+\mu)e_0-(\mu-1)e_{-\mu})\\[5pt]=\!\!\!&0.\end{array}$$
Solving following linear equations on
$e_0,e_{-1},e_{-2},e_{-3},e_{-4}$:
\begin{eqnarray*}
\Bigg\{\begin{array}{lllll}
\Eqb_1(0)=0,&\Eqb_1(2)=0,&\Eqb_1(3)=0,\\[8pt]
\Eqb_2(1)=0,&\Eqb_2(2)=0,&\Eqb_2(3)=0,
\end{array}
\end{eqnarray*}
we obtain $e_0=e_{-1}=e_{-2}=e_{-3}=e_{-4}=0$. Then
$\Eqb_1(\mu)=0$ becomes $(\mu-1)e_{-\mu}=\mu^2e_{-\mu-1}$, which
implies that $e_\mu=0$ for all $\mu\in\Z$. Similarly, applying
$L_{\a,2}=\frac{1}{\a}[L_{0,2},L_{\a,0}]$ with $\a\ne0$ to $v_\mu$,
we obtain
\begin{eqnarray*}
f_{\a,\mu}= \left\{\begin{array}{ll}
\frac{1}{\a}\mu(f_{\mu}-f_{\a+\mu}) &\mbox{if \ }\mu\ne-\a,\\[8pt]
(a+\a)(f_{0}-f_{-\a}) &\mbox{if \ }\mu=-\a.\end{array}\right.
\end{eqnarray*}
Applying $[L_{\mu,2},L_{1,0}]=(\mu+1) L_{\mu+1,2}$ to $v_{-\mu}$
with $\mu\ne0,\pm1$, and $v_{-\mu-1}$ with $\mu\ne0,-1$
respectively, we obtain
$$\begin{array}{ll}\Eqb_3(\mu):=(1-\mu)(f_{1-\mu}-f_1)+\mu(f_{-\mu}-f_1)=0,\\[6pt]
\Eqb_4(\mu):=(a+1)(f_{-\mu-1}-f_{-1})+\mu(a+\mu+1)f_{-\mu-1}-
\mu(a+\mu)f_{-\mu}-\mu f_0=0.\end{array}$$ Solving following linear
equations on $f_0,f_{\pm1},f_{\pm2},f_{\pm3}$:
\begin{eqnarray*}
\Bigg\{\begin{array}{lllll}
\Eqb_3(2)=0,&\Eqb_3(-2)=0,&\Eqb_3(3)=0,\\[8pt]
\Eqb_4(2)=0,&\Eqb_4(-2)=0,&\Eqb_4(1)=0,&\Eqb_4(-3)=0,
\end{array}
\end{eqnarray*}
we obtain $f_0=f_{\pm1}=f_{\pm2}=f_{\pm3}$. Now, rewriting
$\Eqb_3(\mu)=0$ as $\mu f_{-\mu}-(\mu-1)f_{1-\mu}=f_{-1}$, we see
that $f_\mu$ is a constant. Thus Claim 1 holds.
\vskip4pt

{\it Subcase 2.3}: $V\cong A'_{0,1}\oplus \C v_0$ as $\Vir$-modules. \vskip4pt

For
$\mu\ne0$, we claim that $e_\mu=0$ and $f_\mu$ is a constant. We
also claim that $e_0(f_1-f_0)=0$. Then, together with
\eqref{simple-case1} and the trivial actions of $L_{\a,0}$ on $v_0$, it
is easy to check that $\C v_0$ is a quotient module of $V$ (thus Claim 1 holds). Applying
$L_{\a,1}=[L_{1,1},L_{\a-1,0}]$ to $v_\mu$, we obtain
\begin{eqnarray*}
e_{\a,\mu}= \left\{\begin{array}{ll}
(\mu+\a)e_\mu-(\mu+\a-1)e_{\mu+\a-1} &\mbox{if \ }\mu\ne0,-1,\\[8pt]
\a e_0 &\mbox{if \ }\mu=0,\\[8pt]
-\a e_\a &\mbox{if \ }\mu=-1.\end{array}\right.
\end{eqnarray*}
Applying $[L_{\mu,1},L_{1,0}]=\mu L_{\mu+1,1}$ to $v_{-1}$ with
$\mu\ne0,1$, and $v_{1}$ with $\mu\ne-1,-2$ respectively, we obtain
$$\begin{array}{ll}\Eqc_1(\mu):=(\mu+1)e_{\mu+1}-\mu e_\mu=0,\\[6pt]\Eqc_2(\mu):=(\mu+1)e_{\mu+1}-\mu e_\mu+e_1-2e_2=0.\end{array}$$ Solving
following linear equations on $e_{-1},e_1,e_2,e_3$:
\begin{eqnarray*}
\Bigg\{\begin{array}{lllll}
\Eqc_1(2)=0,&\Eqc_1(-1)=0,\\[8pt]
\Eqc_2(2)=0,&\Eqc_2(0)=0,
\end{array}
\end{eqnarray*}
we obtain $e_{-1}=e_1=e_2=e_3=0$. Then $\Eqc_1(\mu)=0$ implies
$e_\mu=0$ for $\mu\ne0$. Similarly, applying
$L_{\a,2}=\frac{1}{\a}[L_{0,2},L_{\a,0}]$ with $\a\ne0$ to $v_\mu$,
we obtain
\begin{eqnarray*}
f_{\a,\mu}= \left\{\begin{array}{ll}
\frac{1}{\a}(\mu+\a)(f_{\mu}-f_{\a+\mu}) &\mbox{if \ }\mu\ne0,\\[8pt]
0 &\mbox{if \ }\mu=0.\end{array}\right.
\end{eqnarray*}
Applying $[L_{1,2},L_{1,0}]=2 L_{2,2}$ to $v_{-1}$ gives
$f_{-1}=f_1$; applying $[L_{-3,2},L_{1,0}]=-2 L_{-2,2}$ to $v_{1}$
gives $f_{-2}=f_2$. Furthermore, applying
$[L_{\mu,2},L_{1,0}]=(\mu+1) L_{\mu+1,2}$ to $v_{-1}$ with
$\mu\ne0,\pm1$, we obtain $\mu f_\mu-(\mu-1)f_{\mu-1}=f_{-1}$, which
implies $f_\mu=f_1$ for $\mu\ne0$. So the first claim holds. At
last, applying $2[L_{1,1},L_{0,2}]=[[L_{1,1},L_{-1,2}],L_{1,0}]$ to $v_0$
gives the second claim.
 \qed\vs{5pt}

To prove Theorem \ref{ISM}(1) and (2) with $q\ne-\frac{1}{2}$, by
\eqref{simple-case1}, \eqref{simple-case2} and Lemma \ref{Th--1}, it
remains to prove \eqref{simple-extend}  (note that \eqref{tri----}
can be regarded as a special case of \eqref{simple-extend}), which
will be done by Lemmas \ref{Th--2}--\ref{Th--4}. Our philosophy  is
the following: First we show that $L_{\a,1}$ acts as zero for
$\a\ne0$. Then for any $(\a,i)\ne(0,-2q)$, we can always choose
infinite many $\b$'s such that $\b(2q+i)-\a(q+1)\ne0$, and so
$L_{\a,i}=\frac{1}{\b(2q+i)-\a(q+1)}[L_{\a-\b,i-1},L_{\b,1}]$ must
act as zero, which gives \eqref{simple-extend}. To prove Theorem
\ref{ISM}(3) and (2) with $q=-\frac{1}{2}$, we should make full use of
the interesting relations \eqref{Bq-relations}.

\begin{lemm}\label{Th--2}
Suppose $q\ne-\frac12,-1$ and we have case \eqref{simple-case2}.
Then \eqref{simple-extend} holds.
\end{lemm}
\ni{\it Proof.~}~Fix $\a_0\ne0$ and assume $L_{\a_0,1}v_\mu=d_\mu
v_{\a_0+\mu}$ for some $d_\mu\in\C$. We want to prove
\begin{equation}\label{La1-acts-as-zero}
d_\mu=0\ \ \mbox{for all}\ \ \mu\in\Z.
\end{equation}
For convenience, we denote $\mub=\mu+a$ and $\mu^-=\mu-a$.
By (\ref{Bq-block}), for $\b,\g\in\Z$, we have
\begin{eqnarray*}
[L_{\g,0},[L_{\b,0},L_{\a_0,1}]] &=&
  (q(\a_0-\b)-\b)(q(\b+\a_0-\g)-\g)L_{\g+\b+\a_0,1}, \\
\left[L_{\g+\b,0},L_{\a_0,1}\right] &=&
  (q(\a_0-\g-\b)-\g-\b)L_{\g+\b+\a_0,1}.
\end{eqnarray*}
Briefly denote  the coefficients of right-hand sides by
$h^{(1)}_{\b,\g}$ and $h^{(2)}_{\b,\g}$ respectively. Applying the
above two equations to $v_\mu$, by \eqref{simple-case2} we obtain
\begin{equation}\label{equ-element}
qh^{(2)}_{\b,\g}\LHS=h^{(1)}_{\b,\g}\RHS,
\end{equation}
where
\begin{eqnarray*}
\LHS &=& (\b+ \a_0 +\mub+b\g)((\a_0+\mub+b\b)d_\mu
     -(\mub+b\b)d_{\b+\mu})\\
   &&-(\mub+b\g)((\a_0+ \mub + \g+  b\b)d_{\g+\mu}
     -(\mub+\g+b\b)d_{\g+\b+\mu}),\\
\RHS &=& (\a_0+\mub+b(\g+\b))d_\mu
  -(\mub+b(\g+\b))d_{\g+\b+\mu}.
\end{eqnarray*}
Now in (\ref{equ-element}), by replacing $(\g,\b,\mub)$ by
$(\g,\g,\b-\g)$, $(\g,-\g,\b)$ and $(-\g,-\g,\b+\g)$ respectively
with $\g\ne0$ and $\b\in a+\Z$, we obtain the following three equations:
\begin{equation} \label{equs-3variables-case1}
\left\{
\begin{array}{rrrrr}
(f^{(1)}_{-\g,\b'}\!-\!f^{(2)}_{\g,\b'})d_{(\b-\g)^-}&\!\!
 +f^{(3)}_{\g,\b}d_{\b^-}&\!\!\!
 +(f^{(1)}_{-\g,\b}+f^{(2)}_{\g,\b})d_{(\b+\g)^-}= 0, \\[6pt]
f^{(4)}_{\g,\b}d_{(\b-\g)^-}&\!\!
 +f^{(5)}_{\g,\b}d_{\b^-}&
 +f^{(4)}_{-\g,\b}d_{(\b+\g)^-}= 0, \\[6pt]
(f^{(1)}_{\g,\b}+f^{(2)}_{-\g,\b})d_{(\b-\g)^-}&\!\!\!\!
 +\!f^{(3)}_{-\g,\b}\!d_{\b^-}&\!\!\!
 +(f^{(1)}_{\g,\b'}\!-\!\!f^{(2)}_{-\g,\b'})d_{(\b+\g)^-}= 0,
\end{array}
\right.
\end{equation}
where $\b'=\b+\a_0$ and
\begin{eqnarray*}
f^{(1)}_{x_1,x_2}&=&q(q\a_0+2(1+q)x_1)(bx_1-x_2)((b-1)x_1-x_2),\\
f^{(2)}_{x_1,x_2}&=&(q\a_0-(1+q)x_1)(q\a_0-x_1)((2b-1)x_1+x_2),\\
f^{(3)}_{x_1,x_2}&=&2q(q\a_0-2(1+q)x_1)(\a_0+bx_1+x_2)((1-b)x_1-x_2),\\
f^{(4)}_{x_1,x_2}&=&q^2\a_0(bx_1-x_2)(\a_0+(b-1)x_1+x_2),\\
f^{(5)}_{x_1,x_2}&=&\a_0((1+3q+2q^2(1+b-b^2))x_1^2+2q^2(\a_0+x_2)x_2).
\end{eqnarray*}
for any $x_1,x_2\in\C$. Regard (\ref{equs-3variables-case1}) as a
system of linear equations on $d_{(\b-\g)^-}$, $d_{\b^-}$,
$d_{(\b+\g)^-}$, and let $\D^{(1)}_{\b,\g}$ denote the determinant
of coefficients, which is a polynomial on $\b$ and $\g$. Observing
that the total degrees on $\b,\g$ of $f^{(1)},f^{(2)},f^{(3)}$ are
$\le3$, and those of $f^{(4)}, f^{(5)}$ are $\le2$. Hence
$\deg\,\D^{(1)}_{\b,\g}\leq8$. Let $P(i,j)$ denote the coefficient
of $\b^i\g^j$ in $\D^{(1)}_{\b,\g}$. Direct computation shows
\begin{eqnarray*}
\D^{(1)}_{\b,\g}&=&P(0,8)\g^8+P(1,6)\b\g^6+P(0,6)\g^6, \mbox{ \ where}\\
P(0,8)&=&8b(1-b)(2b-1)q(1+q)^3(1+2q)\a_0,\\
P(1,6)&=&2(1+q)^2(1+2q)\left(1+q-2q^2+12bq^2-12b^2q^2\right)\a_0^2.
\end{eqnarray*}

If $b\ne\frac{1}{2}$, then $P(0,8)\ne 0$ and thus
$\D^{(1)}_{\b,\g}\ne0$, which implies \eqref{La1-acts-as-zero} holds
by \eqref{equs-3variables-case1}.

If $b=\frac{1}{2}$, then we have
\begin{equation*} \label{p16-b=1/2}
P(1,6)|_{b=\frac{1}{2}}=2(1+q)^2(1+2q)(1+q+q^2)\a_0^2.
\end{equation*}
We use the symbol $\sqrt{-1}$ to stand for the imaginary unit. Then
the primitive cube roots of unity can be written as
$\o=-\frac{1}{2}+\frac{\sqrt{-3}}{2}$,
$\o^2=-\frac{1}{2}-\frac{\sqrt{-3}}{2}$. Suppose $q\ne\o$, $\o^2$,
then $P(1,6)|_{b=\frac{1}{2}}\ne0$, which also gives
\eqref{La1-acts-as-zero}. Now, suppose $q=\o$ or $\o^2$. Canceling
the common term $d_{(\b+\g)^-}$ in the first two equations in
\eqref{equs-3variables-case1} gives an equation (denoted A) on
$d_{(\b-\g)^-}$ and $d_{\b^-}$. In A, replacing $(\b,\g)$ by
$(\b,1), (\b-1,1)$ and $(\b,2)$ respectively gives three equations
(denoted B) on $d_{(\b-2)^-}, d_{(\b-1)^-}$ and $d_{\b^-}$. In B,
canceling the common terms $d_{(\b-2)^-}, d_{(\b-1)^-}$ gives $F_\b
d_{\b^-}=0$, where $F_\b$ is a polynomial on $\b$. Let $H(i)$ denote
the coefficient of $\b^i$ in $F_\b$. In particular, since
$\a_0\in\Z^*$, we have
\begin{eqnarray*}
H(4)|_{q=\o}&=&-3\left(1+\sqrt{-3}\right)
\left(24+16\sqrt{-3}+\left(9-\sqrt{-3}\a_0^2\right)\right)\a_0^3\ne0,\\
H(4)|_{q=\o^2}&=&-3\left(1-\sqrt{-3}\right)
\left(24-16\sqrt{-3}+\left(9+\sqrt{-3}\a_0^2\right)\right)\a_0^3\ne0,
\end{eqnarray*}
each of which again implies \eqref{La1-acts-as-zero} holds.
Hence, by our  philosophy stated before Lemma \ref{Th--2}, \eqref{simple-extend} holds. \qed

\begin{lemm}\label{Th--3}
Suppose $q\ne-\frac12,-1$ and we have case \eqref{simple-case1}.
Then \eqref{simple-extend} holds.
\end{lemm}
\ni{\it Proof.~}~Using a similar argument as
\eqref{equs-3variables-case1}, for any $\b\ne0,\pm\g\,(\g\ne0)$, we
have
\begin{equation}\label{equs-3variables-case2}
\left\{
\begin{array}{rrrrr}
(g^{(1)}_{-\g,\b'}\!-\!g^{(2)}_{\g,\b'})d_{\b-\g}&\!\!
 +g^{(3)}_{\g,\b}d_{\b}&\!\!\!
 +(g^{(1)}_{-\g,\b}+g^{(2)}_{\g,\b})d_{\b+\g}= 0, \\[6pt]
g^{(4)}_{\g,\b}d_{\b-\g}&\!\!
 +g^{(5)}_{\g,\b}d_{\b}&
 +g^{(4)}_{-\g,\b}d_{\b+\g}= 0, \\[6pt]
(g^{(1)}_{\g,\b}+g^{(2)}_{-\g,\b})d_{\b-\g}&\!\!\!\!
 +\!g^{(3)}_{-\g,\b}\!d_{\b}&\!\!\!
 +(g^{(1)}_{\g,\b'}\!-\!\!g^{(2)}_{-\g,\b'})d_{\b+\g}= 0,
\end{array}
\right.
\end{equation}
where $g^{(i)}=f^{(i)}|_{b=1}$ for $1\le i\le5$. Regard
(\ref{equs-3variables-case2}) as a system of linear equations on
$d_{\b-\g}$, $d_{\b}$, $d_{\b+\g}$, and denote $\D^{(2)}_{\b,\g}$
 the determinant of coefficients, $Q(i,j)$ the
coefficient of $\b^i\g^j$ in $\D^{(2)}_{\b,\g}$. Then
\begin{eqnarray*}
\D^{(2)}_{\b,\g}&=&Q(1,6)\g^6+Q(0,6)\g^6 \mbox{ \ where}\\
Q(1,6)&=&2(1-q)(1+q)^2(1+2q)^2\a_0^2,\\
Q(0,6)&=&(1-q)(1+q)^2(1+2q)^2\a_0^3.
\end{eqnarray*}
If $q\ne1$, then both $Q(1,6)$ and $Q(0,6)$ are not equal to zero,
and therefore $\D^{(2)}_{\b,\g}\ne0$ for $\b\ne0,\pm\g$, which
implies $d_\mu=0$ for $\mu\ne0,\pm\g,\pm2\g$ by
\eqref{equs-3variables-case2}. By the arbitrariness of $\g$, we have
$d_\mu=0$ for all $\mu\in\Z^*$. If $q=1$, by applying
$[L_{0,0},L_{\a_0,1}]=\a_0 L_{\a_0,1}$ to $v_\mu$, we obtain $\a_0
d_\mu=0$, which also gives $d_\mu=0$. So \eqref{simple-extend}
always holds by our  philosophy stated before Lemma \ref{Th--2}. \qed

\begin{lemm}\label{Th--4}
Suppose $q=-\frac12$. Then \eqref{simple-extend} holds.
\end{lemm}
\ni{\it Proof.~}~\ \ Recall that
$\BBB(-\frac{1}{4})\raisebox{-2pt}{$\ ^{\sc\ra}_{\sc\ne}\
$}\BBB(-\frac{1}{2})$ in the sense that $\BBB(-\frac{1}{2})$
contains the subalgebra with basis
$\{L'_{\a,i}=\frac{1}{2}L_{\a,2i}\,|\,\a\in\Z, i\in\Z_+\}$
isomorphic to $\BBB(-\frac{1}{4})$. By Lemma \ref{Th--1}, $V$
remains irreducible when regarded as a $\Vir$-module or a
$\BBB(-\frac{1}{4})$-module. By Lemmas \ref{Th--2} and \ref{Th--3},
for any $v_\mu\in V_\mu$, $\a\in\Z$, we have
\begin{equation}\label{Bq-1/2-even}
L_{\a,2i}v_\mu=2L'_{\a,i}v_\mu=0\mbox{\ \ if \ }i\ge1.
\end{equation}
Since $\BBB(-\frac{1}{2})$ can be generated by
$\{L_{0,1},L_{\a,2i}\,|\,\a\in\Z,\,i\in\Z_+\}$, it suffices to
determine the action of $L_{0,1}$. Suppose $L_{0,1}v_\mu=e_\mu
v_\mu$. We claim that $e_\mu$ is a constant (denoted
$s$), which together with \eqref{Bq-1/2-even} gives
\eqref{simple-extend} with $q=-\frac{1}{2}$.

First suppose we have case \eqref{simple-case2} with $q=-\frac{1}{2}$.
For $\a\ne0$, applying $L_{\a,1}=\frac{2}{\a}[L_{0,1},L_{\a,0}]$ to
$v_\mu$, we obtain
$L_{\a,1}v_\mu=\frac{1}{\a}(a+\mu+b\a)(e_\mu-e_{\a+\mu})v_{\a+\mu}$.
Furthermore, by \eqref{Bq-1/2-even} with $i=1$, applying
$L_{\a,2}=\frac{2}{\a}[L_{0,1},L_{\a,1}]$ to $v_\mu$ gives
$\frac{2}{\a^2}(a+\mu+b\a)(e_\mu-e_{\a+\mu})^2=0$. In this equation,
replacing $(\a,\mu)$ by $(1,\mu)$, $(-1,\mu+1)$ and $(2,\mu)$
respectively, we obtain
\begin{eqnarray}
(a+\mu+b)(e_\mu-e_{\mu+1})^2&=&0, \label{B-1/2-case1-1}\\
(a+\mu-b+1)(e_\mu-e_{\mu+1})^2&=&0, \label{B-1/2-case1--1}\\
(a+\mu+2b)(e_\mu-e_{\mu+2})^2&=&0. \label{B-1/2-case1-2}
\end{eqnarray}
If $b\ne\frac{1}{2}$, by comparing \eqref{B-1/2-case1-1} and
\eqref{B-1/2-case1--1}, then we have $(2b-1)(e_\mu-e_{\mu+1})^2=0$,
which gives the claim. If $b=\frac{1}{2}$ and
$a+\frac{1}{2}\notin\Z$, then the claim still holds by
\eqref{B-1/2-case1-1}. If $b=\frac{1}{2}$ and $a+\frac{1}{2}\in\Z$,
denoting $\mu_0=-a-\frac{1}{2}$, then by \eqref{B-1/2-case1-1} we
have, for some $s,s'\in\C$,
\begin{equation}\label{SMSMS}
e_\mu=\Bigg\{
\begin{aligned}
&s&&\mbox{if}\ \ \ \mu\ge\mu_0+1,\\
&s' &&\mbox{if}\ \ \ \mu\le\mu_0.
\end{aligned}
\end{equation}
On the other hand, taking $\mu=\mu_0$ in \eqref{B-1/2-case1-2} gives
$s=s'$, which again gives the claim.

Now suppose we have case \eqref{simple-case1} with $q=-\frac{1}{2}$. By
 similar arguments to those in obtaining \eqref{SMSMS}, for $\mu\ne0,-1$, we
have $(\mu+1)(e_\mu-e_{\mu+1})^2=0$, which gives, for some $s,s'\in\C$,
\begin{equation*}
e_\mu=\Bigg\{
\begin{aligned}
&s &&\mbox{if}\ \ \ \mu\ge1,\\
&s' &&\mbox{if}\ \ \ \mu\le-2.
\end{aligned}
\end{equation*}
Applying $L_{-1,2}=2[L_{-1,1},L_{0,1}]$ to $v_{-2}$ gives
$2(e_{-2}-e_{-1})^2=0$, and so $e_{-1}=s'$.
Furthermore, applying $L_{2,2}=[L_{0,1},L_{2,1}]$ to $v_{-1}$, we
obtain $(e_{-1}-e_1)^2=0$, and so $s=s'$, which gives the claim.\qed\vspb

\ni{\it Proof of Theorem \ref{ISM}(3).}\ \ By \eqref{simple-case1},
\eqref{simple-case2} and Lemma \ref{Th--1}, it remains to prove
\eqref{simple-extend} with $q=-1$, and
\eqref{simple-extend-q=-1-irr}. Similarly to Lemma \ref{Th--4},
recall the relation $\BBB(-\frac{1}{2})\raisebox{-2pt}\
^{\sc\ra}_{\sc\ne}\ \BBB(-1)$ in the sense that $\BBB(-1)$ contains
the subalgebra with basis
$\{L''_{\a,i}=\frac{1}{2}L_{\a,2i}\,|\,\a\in\Z, i\in\Z_+\}$
isomorphic to $\BBB(-\frac{1}{2})$. By Lemma \ref{Th--1}, $V$
remains irreducible when regarded as a $\Vir$-module or a
$\BBB(-\frac{1}{2})$-module. By Lemma \ref{Th--4}, for any
$\a\in\Z$, $i\ge1$, we have, for some $s\in\C$,
\begin{eqnarray}\label{Bq-1-even}
L_{\a,2i}v_\mu = 2L''_{\a,i}v_\mu= \bigg\{\begin{array}{ll}
s v_\mu &\mbox{if \ }(\a,i)=(0,1),\\[8pt]
0 &\mbox{if \ }(\a,i)\ne(0,1)\mbox{ and }i\ge1.\end{array}
\end{eqnarray}
Since $\BBB(-1)$ can be generated by
$\{L_{1,1},L_{\a,2i}\,|\,\a\in\Z,\,i\in\Z_+\}$, it suffices to
determine the action of $L_{1,1}$. Suppose
$L_{\a,1}v_\mu=f_{\a,\mu}v_{\a+\mu}$, and write $f_{1,\mu}=f_\mu$
for short.

First suppose we have case \eqref{simple-case2} with $q=-1$. We claim that
$f_\mu$ is a constant (denoted $t$), which, together with
\eqref{Bq-1-even}, gives \eqref{simple-extend-q=-1-irr}. For any
$\a,\mu\in\Z$, applying $L_{\a,1}=[L_{1,1},L_{\a-1,0}]$ to $v_\mu$
gives
\begin{equation*}\label{La1-vu}
f_{\a,\mu}=(a+\mu+1+b(\a-1))f_\mu-(a+\mu+b(\a-1))f_{\a+\mu-1}.
\end{equation*}
Now applying $[L_{0,1},L_{\a,0}]=0$, $[L_{1,1},L_{\a-1,1}]=0$ to $v_\mu$ gives respectively
$$\begin{array}{lll}
\Eq_1\,(\a,\mu):=(a+\mu+b\a)(f_{0,\mu}-f_{0,\a+\mu})=0;\\[5pt]
\Eq_2\,(\a,\mu):=f_{\a-1,\mu}f_{\a+\mu-1}-f_{\a-1,\mu+1}f_\mu=0;\end{array}$$
applying $[L_{\mu,1},L_{1,0}]=\mu L_{\mu+1,1}$ to $v_\a$ gives
$$\Eq_3\,(\a,\mu):=(a+\mu+\a+b)f_{\mu,\a}-(a+\a+b)f_{\mu,\a+1}-\mu
f_{\mu+1,\a}=0.$$ Solving following equations on
$f_{\mu+1},f_{\mu},f_{\mu-1},f_{\mu-2}$:
\begin{eqnarray}\label{Bq-1-L11}
\left\{\begin{array}{ll}
\Eq_1\,(1,\mu)=0,&\Eq_1\,(1,\mu-1)=0,\\[8pt]
\Eq_2\,(0,\mu)=0,&\Eq_2\,(3,\mu-1)=0,\\[8pt]
\Eq_3\,(0,\mu)=0,&\Eq_3\,(-1,\mu+1)=0,
\end{array}\right.
\end{eqnarray}
we obtain the following possible solutions:
$$\begin{array}{lll}
(\rm i)& f_\mu\mbox{ is a constant for all }\mu\in\Z,\\[8pt]
(\rm ii)& b=0\mbox{ and }f_\mu=\bigg\{\begin{array}{ll}
0 &\mbox{if \ }\mu\ne-a-1,\\[3pt]
t_0 &\mbox{if \ }\mu=-a-1\mbox{ (for some $t_0\ne0$)},\end{array}\\[12pt]
(\rm iii)& b=1\mbox{ and }f_\mu=\bigg\{\begin{array}{ll}
0 &\mbox{if \ }\mu\ne-a,\\[3pt]
t_1 &\mbox{if \ }\mu=-a\mbox{ (for some $t_1\ne0$)},\end{array}
\end{array}$$
Recall that if $V\cong A_{a,b}$ as a $\Vir$-module, then $b=0$ or
$1$ implies $a\notin\Z$. So cases (\rm ii) and (\rm iii) become a
special case of (\rm i), and therefore the claim holds.

Suppose we have case \eqref{simple-case1} with $q=-1$. We claim that
$f_\mu=0$ for $\mu\in\Z^*$, which, together with \eqref{Bq-1-even},
gives \eqref{simple-extend} with $q=-1$. Obviously, $f_{-1}=0$. By
similar arguments to those in obtaining \eqref{Bq-1-L11} (or equivalently taking
$a=0,b=1$ in the last two equations on $f_\mu,f_{\mu-1}$ with
$\mu\neq0,\pm1$ in \eqref{Bq-1-L11}), the claim holds. \qed

\begin{rema}\rm  We can say something more about
the irreducible $\BB$-module of the intermediate series for
$q\in\frac{1}{2}\Z^*_-$.
\baselineskip3pt\lineskip7pt\parskip-3pt
\begin{itemize}\parskip-3pt
\item[{\rm(1)}] Using the relation $\BBB(-\frac{1}{2})\ra\BBB(-\frac{k}{2})$,
 where $k\in\Z_+^*$, one can deduce that $\BBB(-\frac{k}{2})$-module
 $A'_{0,1}(s)$ is also a $\BBB(-\frac{1}{2})$-module
 $A'_{0,1}(\frac{s}{k})$ for any $s\in\C$;
\item[{\rm(2)}] Note that $\BB$ is perfect if and only if
 $q\notin\frac{1}{2}\Z^*_-$. So, by Theorem \ref{ISM},
 an irreducible $\BB$-module of the intermediate series is a trivial extension
 from an irreducible $\Vir$-module of the intermediate series if and only if $\BB$ is perfect.
\end{itemize}
\end{rema}


\small\def\bf{}
\parskip=-1pt\baselineskip=3pt\lineskip=3pt

\end{CJK*}

\begin{thebibliography}{9999}\vskip0pt\small
\parindent=2ex\parskip=-1pt\baselineskip=-1pt
\def\re#1{\bibitem{#1}\label{#1}}


\re{B} R. Block, On torsion-free abelian groups and Lie algebras,
   \textit{Proc. Amer. Math. Soc.} \textbf{9} (1958) 613--620.

\re{DZ} D. Dokovic, K. Zhao, Derivations, isomorphisms and
  second cohomology of generalized Block algebras,
  \textit{Algebra Colloquium} \textbf{3} (1996) 245--272.

\re{KL} V. Kac, J. Liberati, Unitary quasifinite representations of
  $W_\infty$, \textit{Lett. Math. Phys.} \textbf{53} (2000) 11--27.

\re{KR} V. Kac, A. Radul, Quasifinite highest weight modules over
  the Lie algebra of differential operators on the circle,
  \textit{Comm. Math. Phys.} \textbf{157} (1993) 429--457.


\re{LZ} R. Lu, K. Zhao, Classification of irreducible weight modules
  over higher rank Virasoro algebras, \textit{Adv. Math.}
  \textbf{206} (2006) 630-656.

\re{M} O. Mathieu, Classification of Harish-Chandra modules over the
  Virasoro Lie algebra, \textit{Invent. Math.} \textbf{107} (1992) 225--234.

\re{MP} C. Martin, A. Piard, Indecomposable modules over the
  Virasoro Lie algebra and a conjecture of V. Kac,
  \textit{Comm. Math. Phys.} \textbf{137} (1991) 109--132.

\re{OZ} J.M. Osborn, K. Zhao, Infinite-dimensional Lie algebras
  of generalized Block type,  \textit{Proc. Amer. Math. Soc.}
  \textbf{127} (1999) 1641--1650.

\re{S0} Y. Su,
A classification of indecomposable $sl_2(\C)$-modules
and a conjecture of Kac on irreducible modules over the Virasoro algebra,
{\it J. Algebra} {\bf161} (1993), 33--46.

\re{S1} Y. Su, Classification of quasifinite modules over the Lie algebras
  of Weyl type. \textit{Adv. Math.} \textbf{174} (2003) 57--68.

\re{S2} Y. Su, Classification of Harish-Chandra modules over the higher rank
  Virasoro algebras, \textit{Comm. Math. Phys.} \textbf{240} (2003)
  539--551.

\re{S3} Y. Su, Quasifinite representations of a Lie algebra of Block
  type, \textit{J. Algebra} \textbf{276} (2004) 117--128.

\re{S4} Y. Su, Quasifinite representations of a family of Lie algebras of Block type,
  \textit{J.~Pure Appl. Algebra} \textbf{192} (2004) 293--305.

\re{S5} Y. Su, Quasifinite representations of some Lie algebras
  related to the Virasoro algebra, in: Recent developments in algebra and related areas, 213--238,
  Adv. Lect. Math. (ALM), 8, Int. Press, Somerville, MA, 2009.

\re{SX} Y. Su, B. Xin, Classification of quasifinite
  $W_{\infty}$-module, \textit{Israel J. Math.} \textbf{151} (2006) 223--236.

\re{WT} Q. Wang, S. Tan, Quasifinite modules of a Lie algebra
  related to Block type, \textit{J. Pure Appl. Algebra} \textbf{211} (2007)
  596--608.

\re{X1} X. Xu, Generalizations of Block algebras,
  \textit{Manuscripta Math.} \textbf{100} (1999) 489--518.

\re{X2} X. Xu, Quadratic conformal superalgebras,
  \textit{J. Algebra} \textbf{224} (2000) 1--38.

\re{XYZ} C. Xia, T. You, L. Zhou, Structure of a class of Lie
  algebras of Block type, \textit{Comm. Algebra}, in press.

\re{ZZ} H. Zhang, K. Zhao, Representations of the Virasoro-like Lie algebra
  and its $q$-analog, \textit{Comm. Algebra} \textbf{24} (1996) 4361--4372.

\re{Z} K. Zhao, A class of infinite dimensional simple Lie algebras,
  \textit{J. London Math. Soc.} (2) \textbf{62} (2000) 71--84.

\re{ZM} L. Zhu, D. Meng, Structure of degenerate Block algebras,
  \textit{Algebra Colloquium} \textbf{10} (2003) 53--62.


\end{thebibliography}
\end{document}